\documentclass[10pt]{amsart}
\usepackage{amsfonts} 
\usepackage{graphicx,latexsym,bm,amsmath,amssymb,verbatim,multicol,lscape}
\newtheorem{thm}{Theorem} [section]

\theoremstyle{definition}
\newtheorem{defn}[thm]{Definition}
\theoremstyle{remark}

\numberwithin{equation}{section}
\newenvironment{narrow}[2]{%
\begin{list}{}{%
\setlength{\topsep}{0pt}%
\setlength{\leftmargin}{#1}%
\setlength{\rightmargin}{#2}%
\setlength{\listparindent}{\parindent}%
\setlength{\itemindent}{\parindent}%
\setlength{\parsep}{\parskip}}%
\item[]}{\end{list}}

\begin{document}
\title{On deep holes of standard Reed-Solomon codes}
\author{Rongjun Wu}
\address{Mathematical College, Sichuan University, Chengdu 610064, P.R. China}
\email{eugen\_woo@163.com}
\author{Shaofang Hong*}
\address{Yangtze Ctenter of Mathematics, Sichuan University, Chengdu 610064, P.R. China and
Mathematical College, Sichuan University, Chengdu 610064, P.R. China}
\email{sfhong@scu.edu.cn, s-f.hong@tom.com, hongsf02@yahoo.com}
\thanks{*Corresponding author and was supported partially by National Science
Foundation of China Grant \# 10971145 and by the Ph.D. Programs
Foundation of Ministry of Education of China Grant \#20100181110073}

\date{\today}%
\begin{abstract}
Determining deep holes is an important open problem in decoding
Reed-Solomon codes. It is well known that the received word is
trivially a deep hole if the degree of its Lagrange interpolation
polynomial equals the dimension of the Reed-Solomon code. For the
standard Reed-Solomon codes $[p-1, k]_p$ with $p$ a prime, Cheng and
Murray conjectured in 2007 that there is no other deep holes except
the trivial ones. In this paper, we show that this conjecture is not
true. In fact, we find a new class of deep holes for standard
Reed-Solomon codes $[q-1, k]_q$ with $q$ a prime power of $p$. Let
$q \geq 4$ and $2 \leq k\leq q-2$. We show that the received word
$u$ is a deep hole if its Lagrange interpolation polynomial is the
sum of monomial of degree $q-2$ and a polynomial of degree at most
$k-1$. So there are at least $2(q-1)q^k$ deep holes if $k \leq q-3$.
\end{abstract}

\maketitle

\section{\bf Introduction and the statement of the main result}

Let ${\bf F}_q$ be the finite field of $q$ elements with
characteristic $p$. Let $D=\{x_1, ..., x_n\}$ be a subset of ${\bf
F}_q$, which is called the {\it evaluation set}. The {\it
generalized Reed-Solomon code} $\mathcal {C}_q(D,k)$ of length $n$
and dimension $k$ over ${\bf F}_q$ is defined as follows:
$$
\mathcal {C}_q(D,k) = \{(f(x_1), ..., f(x_n)) \in {\bf F}_q^n | f(x)
\in {\bf F}_q[x], {\rm deg} f(x) \leq k - 1\}.
$$
If $D={\bf F}_q^*$, it is called {\it standard Reed-Solomon code}.
We refer the above definition as the {\it polynomial code} version
of the standard Reed-Solomon code. If $D={\bf F}_q$, it is called
{\it extended Reed-Solomon code}. For any $[n,k]_q$ linear code
$\mathcal {C}$, the {\it minimum distance} $d(\mathcal {C})$ is
defined by
$$
d(\mathcal {C}):= {\min}\{d(x,y)|x \in \mathcal {C}, y \in \mathcal
{C}, x \neq y\},
$$
where $d(\cdot,\cdot)$ denotes the {\it Hamming distance} of two
words which is the number of different entries of them and
$w(\cdot)$ denotes the {\it Hamming weight} of a word which is the
number of its non-zero entries. Thus we have
$$
d(\mathcal {C})={\min}\{d(x,0)|0 \neq x \in \mathcal {C}\}
={\min}\{w(x)|0 \neq x \in \mathcal {C}\}.
$$
The {\it error distance} to code $\mathcal {C}$ of a received word
$u \in {\bf F}_q^n$ is defined by $d(u, \mathcal {C}):={\rm
min}\{d(u, v)|v \in \mathcal {C}\}$. Clearly $d(u,\mathcal {C})=0$
if and only if $u \in \mathcal {C}$. The {\it covering radius}
$\rho(\mathcal {C})$ of code $\mathcal {C}$ is defined to be ${\rm
max}\{d(u,\mathcal {C})|u \in {\bf F}_p^n\}$. For generalized
Reed-Solomon code $\mathcal {C}=\mathcal {C}_q(D,k)$, we have that
the minimum distance $d(\mathcal {C})=n-k+1$ and the covering radius
$\rho(\mathcal {C})=n-k$. The most important algorithmic problem in
coding theory is the maximum likelihood decoding (MLD): Given a
received word, find a word $v \in \mathcal {C}$ such that
$d(u,v)=d(u,\mathcal {C})$ \cite{[LW]}. Therefore, it is very
crucial to decide $d(u,\mathcal {C})$ for the word $u$. Sudan
\cite{[Su]} and Guruswami-Sudan \cite{[GS]} provided a polynomial
time list decoding algorithm for the decoding of $u$ when
$d(u,\mathcal {C})\leq n-\sqrt{nk}$. When the error distance
increases, the decoding becomes NP-complete for generalized
Reed-Solomon codes \cite{[GV]}.

When decoding the generalized Reed-Solomon code $\mathcal C$, for a
received word $u=(u_1, ..., u_n)\\ \in {\bf F}_q^n$, we define the
{\it Lagrange interpolation polynomial} $u(x)$ of $u$ by
$$
u(x):=\sum_{i=1}^n u_i \prod_{\substack{j =1\\j \neq
i}}^{n}\frac{x-x_j} {x_i-x_j} \in {\bf F}_q[x],
$$
i.e., $u(x)$ is the unique polynomial of degree at most $n-1$ such
that $u(x_i)=u_i$ for $1 \leq i \leq n$. For $u \in {\bf F}_q^n$, we
define the degree of $u(x)$ to be the {\it degree} of $u$, i.e.,
${\rm deg}(u) ={\rm deg}(u(x))$. It is clear that $d(u,\mathcal
{C})=0$ if and only if ${\rm deg}(u) \leq k-1$. Evidently, we have
the following
simple bounds.\\

\noindent{\bf Lemma 1.1. \cite{[jW]}} {\it For $k \leq {\rm deg}(u)
\leq n-1$, we have the inequality}
$$
n-{\rm deg}(u) \leq d(u,\mathcal
{C}) \leq n-k = \rho.
$$

Let $u \in {\bf F}_q^n$. If $d(u,\mathcal {C})=n-k,$ then the word
$u$ is called a {\it deep hole}. If ${\rm deg}(u)=k$, then the upper
bound is equal to the lower bound, and so $d(u,\mathcal {C})=n-k$
which implies that $u$ is a deep hole. This immediately gives
$(q-1)q^k$ deep holes. We call these deep holes {\it the trivial }
deep holes. It is an interesting open problem to determine all deep
holes. Cheng and Murray \cite{[CM]} showed that for the standard
Reed-Solomon code $[p-1, k]_p$ with $k < p^{1/4-\epsilon}$, the
received vector $(f(\alpha ))_{\alpha\in {\bf F}_p^*}$ cannot be a
deep hole if $f(x)$ is a polynomial of degree $k+d$ for $1\le
d<p^{3/13-\epsilon}$. Based on this result they conjectured that
there is no other deep holes except the trivial ones mentioned
above. Li and Wan \cite{[LW]} use the method of character sums to
obtain a bound on the non-existence of deep holes for extended
Reed-Solomon code $\mathcal {C}_q({\bf F}_q,k)$.

In this paper, we introduce a new method to investigate the standard
Reed-Solomon code $\mathcal {C}_q({\bf F}_q^*,k)$ and particularly
we study its deep holes. In fact, we use the BCH code and the
discrete Fourier transform to find some new deep holes.
We have the following result.\\

\noindent{\bf Theorem 1.2.} {\it Let $q \geq 4$ and $2 \leq k \leq
q-2$. If the Lagrange interpolation polynomial of a received word $u
\in {\bf F}_q^{q-1}$ is of the form either $ax^{q-2}+f_{\le k-1}(x)$
or $ax^k+f_{\le k-1}(x)$, where
$a\in {\bf F}^*_q$ and $f_{\le k-1}(x)\in {\bf F}_q[x]$ is of degree
at most $k-1$, then $u$ is a deep hole of the standard Reed-Solomon
code $\mathcal {C}_q({\bf F}_q^*,k)$.}\\

If one picks $k=q-2$, then the deep holes given by Theorem 1.2 are the same as
the trivial ones. If $k < q-2$, then there exist two types deep
holes, namely, we give $2(q-1)q^k$ deep holes. Theorem 1.2 also
implies that the Cheng-Murray conjecture \cite{[CM]} is not true.

This paper is organized as follows. First we give deep holes for the
cyclic code version of the standard Reed-Solomon codes in Section 2.
Subsequently, in Section 3, by using the results presented in
Section 2 and the discrete Fourier transformation of vectors, we
show Theorem 1.2. Finally, we give in Section 4 some examples of the
received word $u$ which is not a deep hole and whose Lagrange
interpolation polynomial is of the form neither $ax^{q-2}+f_{\le
k-1}(x)$ nor $ax^{k}+f_{\le k-1}(x)$. We suggest a conjecture
on the non-existence of deep holes for the standard Reed-Solomon
code $\mathcal {C}_q({\bf F}_q^*, k)$ at the end of the paper.

\section{\bf Deep hole of the standard Reed-Solomon codes as cyclic codes}

In the present section, we consider the deep hole for the cyclic
code version of the standard Reed-Solomon codes. Throughout this
paper, we let $\alpha$ be a fixed primitive element of ${\bf F}_q$
and
$$
g(x)=(x-\alpha)(x-\alpha^2)...(x-\alpha^{d-1}).
$$
We define the standard Reed-Solomon code with length $n=q-1$ as
follows:
$$
\mathcal{C}:=\{g(x)m(x) \in {\bf F}_q[x] \mid m(x)\in {\bf F}_q[x]
\mbox{ and }\deg(m(x)) \leq k-1 \}. \eqno (2.1)
$$
and refer it as the cyclic code version of the standard Reed-Solomon
code. It is a maximal distance separable code, i.e., we have
$d=n-k+1=q-k$.

Now let $D={\bf F}_q^*=\{1,\alpha, ..., \alpha^{q-2}\}$. Then the
the polynomial code version of the standard Reed-Solomon code can be
stated as follows
$$
\mathcal {C}_q({\bf F}_q^*,k) := \{(f(1), f(\alpha), ...,
f(\alpha^{q-2})) \in {\bf F}_q^n | f(x) \in {\bf F}_q[x], {\rm deg}
f(x) \leq k-1\}.
$$

Let's recall the definition of BCH code and the BCH bound.

\defn{\cite{[vL]} Let $\delta \geq 2$ be an integer.
A cyclic code of length $n$ over ${\bf F}_q$ is called a {\it BCH
code of designed distance $\delta$} if its generator $g(x)$ is the
least common multiple of the minimal polynomials of $\beta^{l},
\beta^{l+1}, ..., \beta^{l+\delta-2}$ for some positive integer $l$,
where $\beta$ is a primitive $n$th root of unity. If $n=q^m - 1$ for
some integer $m \geq 1$, then the BCH code is called {\it primitive}. }\\
\\
The designed distance is explained by the following result which is
usually called {\it BCH bound}.\\

\noindent{\bf Lemma 2.2.} {\it \cite{[vL]} The minimum distance of a
BCH code with designed distance $\delta $ is at least $\delta$.}\\

We are now in a position to give the main result of this section.\\

\noindent{\bf Theorem 2.3.} {\it For the cyclic code version of the
standard Reed-Solomon code $\mathcal{C}$ defined as in (2.1), let
$g_1(x):=g(x)/(x-\alpha)$. Then $a g_1(x)+l(x)g(x)$ is a deep hole
of $\mathcal{C}$ for any $a \in {\bf F}_q^*$ and $l(x)\in {\bf
F}_q[x]$ with ${\rm deg}(l(x))\leq k-1$.}\\
\\
{\it Proof.} First we show that $g_1(x)$ is a deep
hole. For any codeword $m(x)g(x)\in \mathcal {C}$ with ${\rm
deg}(m(x)) \leq k-1$, we have
$$
d(g_1(x),m(x)g(x))=w(g_1(x)(1-m(x)(x-\alpha)))
$$
where $w(\cdot)$ means the number of the nonzero coefficients of the
polynomial. Let
$$
\Delta(x):=1-m(x)(x-\alpha).
$$
We define the primitive BCH code $C'$ as follows:
$$
\mathcal {C}'=\{ c(x)=\sum_{i=0}^{q-2}c_ix^i \in {\bf
F}_q[x]/(x^{q-1}-1)\big{|}c(\alpha^2)=...=c(\alpha^{d-1})=0 \}.
$$
Then the designed distance $\delta$ of $\mathcal{C'}$ is $d-1$ since
its generator, the least common multiple of the minimal polynomials
of $\alpha^2, ..., \alpha^{d-1}$, is  $g_1(x)$. For the BCH code, by
Lemma 2.2 we know that the minimum distance $d'$ is at least
$\delta$. Since ${\rm deg}(\Delta(x)) \leq k$ and $d=q-k$, we have
${\rm deg}(g_1(x)\Delta(x)) \leq q-2$. But $g_1(\alpha ^i)=0$ for
$2\le i\le d-1$. So $g_1(x)(\Delta(x)) \in \mathcal {C}'$. Note that
$\mathcal {C}'$ holds minimum distance $d'\geq d-1$. Therefore
$w(g_1(x)\Delta(x))\geq d-1=q-1-k$. It then follows that $g_1(x)$ is
a deep hole.

Subsequently, we compute the distance between $a g_1(x)+l(x)g(x)$
and the standard Reed-Solomon codes $\mathcal{C}$. Since $a \in {\bf
F}_q^{*}$ and $l(x) \in {\bf F}_q[x]$, we infer that \begin{align*}
&d(ag_1(x)+l(x)g(x),\mathcal{C})\\
&= {\rm min}\{d(ag_1(x)+l(x)g(x),m(x)g(x))\mid m(x) \in {\bf F}_q[x]
\mbox{ and
} {\rm deg}(m(x))\leq k-1\}\\
&= {\rm min}\{d(ag_1(x)+l(x)g(x),(m(x)+l(x))g(x))\mid m(x) \in {\bf
F}_q[x] \mbox{ and
} {\rm deg}(m(x))\leq k-1\}\\
&= {\rm min}\{d(ag_1(x),m(x)g(x))\mid m(x) \in {\bf F}_q[x] \mbox{
and
} {\rm deg}(m(x))\leq k-1\}\\
&= {\rm min}\{d(ag_1(x),am(x)g(x))\mid m(x) \in {\bf F}_q[x] \mbox{
and
} {\rm deg}(m(x))\leq k-1\}\\
&= {\rm min}\{d(g_1(x),m(x)g(x))\mid m(x) \in {\bf F}_q[x] \mbox{
and
} {\rm deg}(m(x))\leq k-1\}\\
&= d(g_1(x),\mathcal {C})=q-1-k=n-k.
\end{align*}
Therefore $a g_1(x)+l(x)g(x)$ is a deep hole. This completes the
proof of Theorem 2.3.\\

By the similar argument as in the proof of Theorem 2.3, one can
deduce the following result.\\
\\
{\bf Theorem 2.4.} {\it Let $g_2(x)=g(x)/(x-\alpha^{d-1})$. Then $a
g_2(x)+l(x)g(x)$ is a deep hole of $\mathcal{C}$ for any $a \in
{\bf F}_q^*$ and $l(x) \in {\bf F}_q[x]$ with ${\rm deg}(l(x))\leq k-1$.}\\

It is natural to ask: Is there any other deep holes except those in
Theorems 2.3 and 2.4? We believe the answer to this question should
be negative if the characteristic of ${\bf F}_q$ is an odd prime,
but we have not yet found a proof. In the next section,
by discrete Fourier transformation, and Theorem 2.3 we will arrive
at the deep holes of the standard Reed-Solomon code for the
polynomial code version.

\section{\bf Proof of Theorem 1.2}

In this section, we first recall the definitions of the discrete
Fourier transform (DFT) and inverse discrete Fourier transform
(IDFT) of vectors over ${\bf F}_q$ (see, for example, \cite{[Mc]}).
Then we use them to give the relationship between the polynomial
code version and cyclic code version of the standard Reed-Solomon
codes. Let $\alpha$ be the fixed primitive element of ${\bf F}_q$.
Let $V=(V_0, V_1, ..., V_{q-2}) \in {\bf F}_q^{q-1}$. The {\it
discrete Fourier transform (DFT)} of $V$, denoted by
$\hat{V}=(\hat{V}_0, \hat{V}_1, ..., \hat{V}_{q-2})$, is defined as
follows:
$$
\hat{V}_j:=\sum_{i=0}^{q-2}V_i \alpha^{ij} \eqno (3.1)
$$
for $j=0,1, ..., q-2.$ We also call $V$ the {\it inverse discrete
Fourier transform (IDFT)} of $\hat V$. Sometimes we call the $V_i$'s
the "time-domain" coordinates of the vector $V$, and the
$\hat{V}_j$'s the "frequency-domain" coordinates, of the vector $V$.
The time-domain components can be recovered from the
frequency-domain components via the IDFT:
$$
V_i=\frac{1}{q-1}\sum_{j=0}^{q-2}\hat{V}_j \alpha^{-ij} \eqno (3.2)
$$
for $i=0,1, ..., q-2.$  If we interpret the components of the $V$
and $\hat{V}$ as the coefficients of polynomials, i.e., if we define
generating functions $V(x)$ and $\hat{V}(x)$ by $$V(x)=V_0 + V_1x +
\cdots + V_{q-2}x^{q-2}$$ and $$\hat{V}(x)=\hat{V}_0 + \hat{V}_1x +
\cdots + \hat{V}_{q-2}x^{q-2},$$ then the DFT and IDFT relationships
(3.1) and (3.2) can be rewritten in the following:
$$
\hat{V}_j=V(\alpha^j) \eqno (3.3)
$$
and
$$
V_i= \frac{1}{q-1}\hat{V}(\alpha^{-i}). \eqno (3.4)
$$
We call $\hat{V}(x)$ the {\it DFT of $V(x)$} and $V(x)$ the {\it
IDFT of $\hat{V}(x)$}. Evidently, the fact $\hat{V}(x)$ is the DFT
of $V(x)$ implies that $V(x)$ is the IDFT of $\hat{V}(x)$, and vise
versa. Note that if $\widehat{V^{(1)}}$ and $\widehat{V^{(2)}}$ are
the DFT of $(q-1)$-dimensional vectors $V^{(1)}$ and $V^{(2)}$ over
${\bf F}_q$ respectively, then $\lambda \widehat{V^{(1)}} + \mu
\widehat{V^{(2)}}$ is the DFT of $\lambda V^{(1)} + \mu V^{(2)}$ for
any $\lambda \in {\bf F}_q$ and $\mu \in {\bf F}_q$. Note also that
if $V^{(1)}$ and $V^{(2)}$ are the IDFT of $(q-1)$-dimensional
vectors $\widehat{V^{(1)}}$ and $\widehat{V^{(2)}}$ over ${\bf F}_q$
respectively, then $\lambda V^{(1)} + \mu V^{(2)}$ is the IDFT of
$\lambda\widehat{V^{(1)}} + \mu \widehat{V^{(2)}}$ for any $\lambda
\in {\bf F}_q$ and $\mu \in {\bf F}_q$. In other words, the DFT and
the IDFT are linear operations in ${\bf F}_q^{q-1}$.

There is a one-to-one correspondence of codewords between the cyclic
code version and polynomial code version of the standard Reed-Solomon
code as the following lemma shows.\\

\noindent{\bf Lemma 3.1.} {\it {\rm (1).} For any $s(x)\in {\bf
F}_q[x]$ with ${\rm deg}(s(x)) \leq k-1$, there exists an $l(x) \in
{\bf F}_q[x]$ with ${\rm deg}(l(x)) \leq k-1$ such that $l(x)g(x)$
is the DFT of $s(x)$.

{\rm (2).} For any $l(x)\in {\bf F}_q[x]$ with ${\rm deg}(l(x)) \leq
k-1$, there exists an $s(x) \in {\bf F}_q[x]$ with ${\rm deg}(s(x))
\leq k-1$ such that $s(x)$ is the IDFT of $l(x)g(x)$.}

\begin{proof}
(1). Write $s(x)=\sum_{i=0}^{k-1}s_i x^i$ and let
$\hat{s}(x)=\sum_{j=0}^{q-2}\hat{s_j} x^j$ be the DFT of $s(x)$. By
(3.3) we have
$$
\hat{s_j}=s(\alpha^j)
$$
for $0\le j\le q-2$. Thus
$$
\hat{s}(x)=\sum_{j=0}^{q-2}s(\alpha^j) x^j. \eqno (3.5)
$$

Since $g(x)=(x-\alpha)(x-\alpha^2) \cdots (x-\alpha^{d-1})$,
$\deg{\hat{s}(x)} \leq q-2$ and $k-1=q-2-(d-1)$, in order to prove
that there is an $l(x) \in {\bf F}_q[x]$ with ${\rm deg}(l(x)) \leq
k-1$ such that $\hat{s}(x)=l(x)g(x)$, it is sufficient to show that
if $1 \leq m \leq d-1$, then $\hat{s}(\alpha^m)=0$. By (3.5) we have
\begin{align*}
\hat{s}(\alpha^m) &= \sum_{j=0}^{q-2} s(\alpha^j)\alpha^{mj}\\
                  &= \sum_{j=0}^{q-2} \sum_{i=0}^{k-1}s_i
                  \alpha^{ji}\alpha^{mj} \\
                  &= \sum_{i=0}^{k-1}s_i \sum_{j=0}^{q-2}
                  \alpha^{(m+i)j}. \ \ \ \ \ \ \ \ \ \ \ \ \ \ \ \ \ \ \        (3.6)
\end{align*}

Since $d=q-k$, one infers that $1 \leq m+i\le d+k-2=q-2$ for $0 \leq
i \leq k-1$ and $1\le m\le d-1$. This implies that $\alpha^{m+i}
\neq 1$. Since $\alpha$ is the primitive root of ${\bf F}_q$, it
follows immediately that $\sum_{j=0}^{q-2} \alpha^{(m+i)j}=0$. Hence
by (3.6), we get that $\hat{s}(\alpha^m) = 0$ for $1 \leq m \leq
d-1$ as required. Part (1) is proved.

(2). Let $s(x)=\sum_{i=0}^{q-2} s_i x^i$ be the IDFT of $l(x)g(x)$.
Then by (3.4) we have
 $$s_i = \frac{1}{q-1} l(\alpha^{-i})
g(\alpha^{-i}) = \frac{1}{q-1} l(\alpha^{-i}) g(\alpha^{q-1-i})
$$
for $0 \leq i \leq q-2$. But $g(x) = \prod_{i=k}^{q-2}
(x-\alpha^{q-1-i})$. So we have $g(\alpha^{q-1-i}) = 0$ for $k \leq
i \leq q-2$. One then derives that $s_i = 0$ for $k \leq i \leq
q-2$. That is, ${\rm deg}(s(x)) \leq k-1$ as desired. Part (2) is
proved.

The proof of Lemma 3.1 is complete.
\end{proof}

Consequently we describe the DFT of $u(x)=\sum_{i=0}^{q-2}u_{i}
x^{i} \in {\bf F}_q[x]$ with
 $u_{q-2} \neq 0$ and $u_i=0$ for all $k \leq i \leq q-3$ in
the following lemma.\\

\noindent{\bf Lemma 3.2.} {\it Let $u(x)=\sum_{i=0}^{q-2}u_{i} x^{i}
\in {\bf F}_q[x]$ satisfy $u_{q-2} \neq 0$ and $u_i=0$ for all $k
\leq i \leq q-3$. Then there exist $a \in {\bf F}_q^*$ and $l(x) \in
{\bf F}_q[x]$ with degree $\leq k-1$ such that $a g_1(x)+l(x)g(x)$
is the DFT of $u(x)$.}

\begin{proof}
Let
$$
\hat{V}(x)=g_1(x)=(x-\alpha^{q-1-(q-3)})(x-\alpha^{q-1-(q-4)})
\cdots (x-\alpha^{q-1-k})
$$
and write its IDFT as
$$
V(x)=V_0 + V_1x + \cdots + V_{q-2}x^{q-2}.
$$
Then
$$V_i=\frac{1}{q-1}\hat{V}(\alpha^{q-1-i})=\frac{1}{q-1}g_1(\alpha^{q-1-i})$$
for $i=0,1,..., q-2.$ For $i=k, ..., q-3$, since $\alpha^{q-1-i}$ is
a zero of $g_1(x)$, we have $V_i=0$. But
$$V_{q-2}=\frac{1}{q-1}g_1(\alpha)=\frac{1}{q-1}\prod_{j=k}^{q-3}(\alpha-\alpha^{q-1-j})
\neq 0.$$ Thus the IDFT of $g_1(x)$ is
$$
V(x)=V_{q-2}x^{q-2} +
V_{k-1}x^{k-1} + V_{k-2}x^{k-2} + \cdots +V_1x+V_0,
$$
where $V_{q-2}\neq 0$.

On the other hand, by Lemma 3.1 (1) there exist $l_1(x) \in {\bf
F}_q[x]$ and $l_2(x) \in {\bf F}_q[x]$ with $\deg(l_1(x)) \leq k-1$
and $\deg(l_2(x)) \leq k-1$ such that $ l_1(x)g(x)$ equals the DFT
of $\sum_{i=0}^{k-1}V_i x^i $ and $l_2(x)g(x)$ is equal to the DFT
of $\sum_{i=0}^{k-1}u_i x^i$. Since the DFT and IDFT are linear
operations in ${\bf F}_q^{q-1}$, it follows that the IDFT of $g_1(x)
- l_1(x)g(x)$ is $V_{q-2} x^{q-2}$.

Now let $a=\frac{u_{q-2}}{V_{q-2}} \in {\bf F}_q^{*}$ and $l(x) = -
a l_1(x) + l_2(x)$. Then $\deg{l(x)} \leq k-1$. One can deduce that
$$
\mbox{ the IDFT of } a g_1(x)- a l_1(x)g(x) = u_{q-2} x^{q-2}.
$$
Hence
$$
a g_1(x)+l(x)g(x) = \mbox{  the DFT of }u_{q-2} x^{q-2} +
\sum_{i=0}^{k-1}u_i x^i
$$
as required. Lemma 3.2 is proved.
\end{proof}

\noindent{\bf Remark.} Similarly, one can prove that there exist an
$a \in {\bf F}_q^*$ and an $l(x) \in {\bf F}_q[x]$ with degree $\leq
k-1$
such that $ag_2(x)+l(x)g(x)$ is the DFT of $x^k$.\\

Next, we show that the DFT is distance preserved.\\

\noindent{\bf Lemma 3.3.} {\it Let $s(x), t(x) \in {\bf F}_q[x]$ be
of degree at most $q-2$, $s=(s(1),s(\alpha),..., s(\alpha^{q-2}))$
and $t=(t(1), t(\alpha), ..., t(\alpha^{q-2}))$. Let $\hat{s}(x)$
and $\hat{t}(x)$ be the DFT of $s(x)$ and $t(x)$ respectively. Then
we have}
 $$d(s,t)=d(\hat{s}(x),\hat{t}(x)).$$
\begin{proof}
Since $d(s,t) = w(s-t)$ and $d(\hat{s}(x),\hat{t}(x)) = w(\hat{s}(x)
- \hat{t}(x))$, in order to prove that
$d(s,t)=d(\hat{s}(x),\hat{t}(x))$, it suffices to show that
$$w(s-t) = w(\hat{s}(x) - \hat{t}(x)).$$

Note that $s-t = (s(1)-t(1), s(\alpha)-t(\alpha), ...,
s(\alpha^{q-2})-t(\alpha^{q-2}))$. So one can easily check that
$$
w(s-t) = q - 1 - \#\{ \mbox{ ditinct roots in } {\bf F}_q^* \mbox{
of } s(x)-t(x)=0\}.
$$
Thus one needs only to show that
$$
w(\hat{s}(x) - \hat{t}(x)) = q - 1 - \#\{ \mbox{ ditinct roots in }
{\bf F}_q^* \mbox{ of } s(x)-t(x)=0\}. \eqno (3.7)
$$
Let $\hat{v}(x) = \hat{v}_0+\hat{v}_1 x+... +\hat{v}_{q-2}x^{q-2} =
\hat{s}(x) - \hat{t}(x) \in {\bf F}_q[x]$. Since $\hat{s}(x) -
\hat{t}(x) = \widehat{s-t}(x)$, it follows that $v(x) = s(x) - t(x)$
is the IDFT of $\hat{v}(x)$. Therefore $\hat{v}_j = v(\alpha^j)$.
Thus to prove (3.7) is equivalent to show the following identity:
$$
w(\hat{v}(x)) = q-1-\#\{ \mbox{ distinct roots in } {\bf F}_q^*
\mbox{ of  } v(x)=0\}. \eqno (3.8)
$$

But we have
\begin{align*}
w(\hat{v}(x)) &= \#\{ 0 \leq j \leq q-2 \mid \hat{v}_j \neq 0 \}\\
              &= \#\{ 0 \leq j \leq q-2 \mid v(\alpha_j) \neq 0 \}\\
              &= q - 1 - \#\{ 0 \leq j \leq q-2 \mid v(\alpha_j) = 0 \}\\
              &= q - 1 - \#\{ \mbox{ distinct roots in } {\bf F}_q^*
\mbox{ of  } v(x)=0\}.
\end{align*}
So (3.8) is proved. This completes the proof of Lemma 3.3.
\end{proof}

Finally, we show Theorem 1.2 as follows. \\
\\
{\it Proof of Theorem 1.2:} Write $u(x)=\sum_{i=0}^{q-2}u_{i} x^{i}
\in {\bf F}_q[x]$ as the Lagrange interpolation polynomial of the
received word $u$, where either $u_{q-2}\ne 0$ and $u_i = 0$ for $k \leq i\leq q-3$
or $u_{k}\ne 0$ and $u_i = 0$ for $k+1\leq i \le q-2$. Further, we can write $u = (u(1), u(\alpha), ...,
u(\alpha^{q-2}))\in {\bf F}_q^{q-1}$, where $\alpha $ is the fixed
primitive element of ${\bf F}_q$. By Lemma 3.3, we have
\begin{align*}
d(u,\mathcal{C}_q({\bf F}_q^*,k))
&=\min\{d(u,v) \mid v \in
\mathcal{C}_q({\bf F}_q^*,k)\}\\
&=\min\{d(\hat{u}(x),\hat{v}(x)) \mid v \in \mathcal{C}_q({\bf
F}_q^*,k)\}.
\end{align*}
From $v = (v(1), v(\alpha), ..., v(\alpha^{q-2})) \in
\mathcal{C}_q({\bf F}_q^*,k)$, we get that $\deg(v(x)) \leq k - 1$.
Let $\mathcal C$ be defined as in (2.1). Then by Lemma 3.1, we have
$\hat{v}(x)= l(x)g(x) \in \mathcal{C}$. So
$$
d(u,\mathcal{C}_q({\bf F}_q^*,k))=\min\{d(\hat{u}(x),\hat{v}(x))
\mid \hat{v}(x) \in \mathcal{C}\}.
$$

On the other hand, by Lemma 3.2 and the remark after Lemma 3.2, we obtain that
either $\hat{u}(x)= ag_1(x)+l(x)g(x)$ or $\hat{u}(x)= ag_2(x)+l(x)g(x)$
for some $a \in {\bf F}_q^*$ and $l(x) \in {\bf
F}_q[x]$ with degree $\leq k-1$. Then we have either
$$
\min\{d(\hat{u}(x),\hat{v}(x)) \mid \hat{v}(x) \in
\mathcal{C}\}=\min\{d(ag_1(x)+l(x)g(x),\hat{v}(x)) \mid \hat{v}(x)
\in \mathcal{C}\}
$$
or
$$
\min\{d(\hat{u}(x),\hat{v}(x)) \mid \hat{v}(x) \in
\mathcal{C}\}=\min\{d(ag_2(x)+l(x)g(x),\hat{v}(x)) \mid \hat{v}(x)
\in \mathcal{C}\}.
$$
By Theorems 2.3 and 2.4, we know that $ag_1(x)+l(x)g(x)$ and $ag_2(x)+l(x)g(x)$ are both deep holes of $\mathcal{C}$. Then one can deduce that
$$
\min\{d(\hat{u}(x),\hat{v}(x)) \mid \hat{v}(x) \in
\mathcal{C}\}=q-1-k.
$$
Therefore we arrive at $d(u,\mathcal{C}_q({\bf F}_q^*,k))=q-1-k$. Thus
$u$ is a deep hole. This complete the proof of Theorem 1.2. \hfill$\Box$

\section{\bf Examples and conjecture}

Let $q \geq 4$ and $2 \leq k \leq q-2$. For any received word $u \in
{\bf F}_q^{q-1}\setminus\mathcal {C}_q({\bf F}_q^*,k)$, it is clear
that the Lagrange interpolation polynomial of $u$ is of degree no
less than $k$ and no more than $q-2$. By Theorem 1.2 we know that if
the Lagrange interpolation polynomial of the received word $u$ is of
the form either $ax^{q-2}+f_{\le k-1}(x)$ or $ax^k+f_{\le k-1}(x)$,
where $a\in {\bf F}^*_q$ and $\deg f_{\le k-1}(x)\le k-1$, then $u$
is a deep hole of the standard Reed-Solomon code $\mathcal
{C}_q({\bf F}_q^*,k)$. The following two examples tell us that there
exist some received words $u$ which are not deep holes and whose
Lagrange interpolation polynomials are of the form neither
$ax^{q-2}+f_{\le k-1}(x)$ nor $ax^{k}+f_{\le k-1}(x)$. In what
follows, we let $q=11, n=q-1=10, \alpha=2, d=6, k=5$ and
$x_i=2^{i-1}$ for $1\le i\le 10$. So we have
$$
\mathcal{C}_{11}({\bf F}_{11}^*,5) = \{(f(1),f(2), ..., f(2^9)) \in
{\bf F}_{11}^{10} | f(x) \in {\bf F}_{11}[x], {\rm deg} f(x) \leq
4\}.
$$
\\
{\bf Example 4.1.} With received words $u$ as in Table \ref{ex1}
whose Lagrange interpolation polynomial $u(x)=\sum_{i=0}^{q-2}u_i
x^i$ satisfies that $u_{i_0} \neq 0$ for exactly one $k+1 \le i_0
\le q-3$ and $u_i=0$ for all $k \le i\le q-2, i \neq i_0$. Using
Matlab R2009b, we search and find a codeword $v$ as in Table
\ref{ex1} such that $d(u,v)<n-k=5$. But $d(u, \mathcal{C}_{11}({\bf
F}_{11}^*,5))\le d(u,v)$. Therefore $d(u, \mathcal{C}_{11}({\bf
F}_{11}^*, 5))<n-k=5$. Namely, the three received words $u$ in Table
\ref{ex1} are not deep holes.
\begin{table}
\centering
\caption{The received words for Example 4.1}\label{ex1}
\begin{tabular}{|c|c|c|c|}
  \hline
  Received word $u$ & Lagrange interpolation polynomial of $u$ &
  Codeword $v$  & $d(u,v)$\\
\hline\hline
  $(8,8,7,8,1,0,0,0,0,0)$ & $1+x+4x^2+3x^3+6x^4+4x^8$ &
  $(0,8,7,8,8,3,0,0,0,7)$ & 4\\
  \hline
  $(4, 9, 5, 1, 1, 0, 0, 0, 0, 0)$ & $2+10x+3x^2+2x^3+8x^4+x^7$ &
  $(0, 3, 9, 9, 10, 0, 8, 4, 0,  0)$ & 4\\
  \hline
  $(2, 3, 9, 9, 1, 0, 0, 0, 0, 0)$ & $9+2x+2x^2+10x^3+2x^4+10x^6$ &
  $(0, 3, 9, 9, 10, 0, 8, 4, 0, 0)$ & 4\\
  \hline
\end{tabular}
\end{table}\\
\\
\noindent {\bf Example 4.2.} (1). With received words $u$ as in
Table \ref{ex2} whose Lagrange interpolation polynomial
$u(x)=\sum_{i=0}^{q-2}u_i x^i$ satisfies that $u_{q-2} \neq 0$ and
$u_{i_0} \neq 0$ for $k \le i_0 \le q-3$. From Table \ref{ex2}, one
can read that $d(u,\mathcal{C}_{11}({\bf F}_{11}^*,5)) \le
d(u,0)=w(u) \le 4< n-k=5$. In other words, the four received words
$u$ in Table \ref{ex2} are not deep holes.
\begin{table}
\centering
\caption{The received words for Example 4.2
(1)}\label{ex2}
\begin{tabular}{|c|c|c|c|}
  \hline
  Received word $u$ & Lagrange interpolation polynomial of
  $u$ & Weight $w(u)$\\
\hline\hline
  $(8, 1, 2, 9, 0, 0, 0, 0, 0, 0)$ & $2+5x+3x^2+x^3+ x^4+7x^8+9x^9$ & 4\\
  \hline
  $(1, 0, 4, 5, 0, 0, 0, 0, 0, 0)$ & $1+7x+4x^2+2x^3+1x^4-x^7+9x^9$ & 3\\
  \hline
  $(3, 0, 6, 0, 8, 0, 0, 0, 0, 0)$ & $5+4x+5x^2+5x^3+7x^4+x^6+9x^9$ & 3\\
  \hline
  $(4, 1, 4, 3, 0, 0, 0, 0, 0, 0)$ & $-1+x+8x^2+7x^3+6x^4+7x^5+9x^9$ & 4\\
  \hline
\end{tabular}
\end{table}\\
\noindent(2). With received words $u$ as in Table \ref{ex3} whose
Lagrange interpolation polynomial $u(x)=\sum_{i=0}^{q-2}u_i x^i$
satisfies that $u_{q-2} \ne 0$ , $u_{i_1} \ne 0$ and  $u_{i_2} \ne
0$ for $k \le i_1 \ne i_2 \le q-3$. Using Matlab R2009b, we search
and find a codeword $v$ as in Table \ref{ex3} such that
$d(u,v)<n-k=5$. Therefore $d(u, \mathcal{C}_{11}({\bf
F}_{11}^*,5))\le d(u, v)<n-k=5$. That is, the six received words $u$
in Table \ref{ex3} are not deep holes.
\begin{table}
\centering \caption{The received words for Example 4.2
(2)}\label{ex3}
\begin{narrow}{-1in}{-1in}
\begin{tabular}{|c|c|c|c|}
  \hline
  Received word $u$ & Lagrange interpolation polynomial of $u$ &
  Codeword $v$  & $d(u,v)$\\
\hline\hline
  $(6,4,10,4,3,0,0,0,0,0)$ & $6+6x+3x^2+9x^3+ x^4+x^7+4x^8+9x^9$ &
  $(0,4,10,7,3,0,3,0,0,5)$ & 4\\
  \hline
  $(4,9,3,1,3,0,0,0,0,0)$ & $2+9x+2x^2+6x^3+6x^4+10x^6+4x^8+9x^9$ &
  $(2,9,3,7,10,0,0,0,0,10)$ & 4\\
  \hline
  $(3,3,10,6,3,0,0,0,0,0)$ & $8+1x+10x^2+4x^3+7x^4+4x^5+4x^8+9x^9$ &
  $(1,3,10,4,3,0,0,7,0,0)$ & 4\\
  \hline
  $(0,10,1,5,3,0,0,0,0,0)$ & $3+7x+x^2+5x^3+8x^4+10x^6+x^7+9x^9$ &
  $(0,0,0,0,0,0,0,0,0,0)$ & 4\\
  \hline
  $(10,4,8,10,3,0,0,0,0,0)$ & $9+10x+9x^2+3x^3+9x^4+4x^5+x^7+9x^9$ &
  $(1,4,8,10,3,1,0,1,4,0)$ & 4\\
  \hline
  $(8,9,1,7,3,0,0,0,0,0)$ & $5+2x+8x^2+3x^4+4x^5+10x^6+9x^9$ &
  $(2,9,3,7,10,0,0,0,0,10)$ & 4\\
  \hline
\end{tabular}
\end{narrow}
\end{table}\\
(3). With received words $u$ as in Table \ref{ex4} whose Lagrange
interpolation polynomial $u(x)=\sum_{i=0}^{q-2}u_i x^i$ satisfies
that $u_{q-2} \ne 0$, $u_{i_1} \ne 0$, $u_{i_2} \ne 0$ and
$u_{i_3}\ne 0$ for $k \le i_1<i_2<i_3\le q-3$. Using Matlab R2009b,
we search and find a codeword $v$ as in Table \ref{ex4} such that
$d(u,v)<n-k=5$. It follows that $d(u, \mathcal{C}_{11}({\bf
F}_{11}^*,5))\le d(u,v)<n-k=5$, i.e., the four received words $u$ in
Table \ref{ex4} are not deep holes.
\begin{table}
\centering \caption{The received words for Example 4.2
(3)}\label{ex4}
\begin{narrow}{-1in}{-1in}
\begin{tabular}{|c|c|c|c|}
  \hline
  Received word $u$ & Lagrange interpolation polynomial of $u$ &
  Codeword $v$  & $d(u,v)$\\
\hline\hline
  $(8,7,8,2,4,0,0,0,0,0)$ & $4+8x+5x^2+8x^3+3x^4+10x^6+x^7+4x^8+9x^9$ &
  $(1,6,8,2,2,0,8,0,0,0)$ & 4\\
  \hline
  $(7,1,4,7,4,0,0,0,0,0)$ & $10+2x^2+6x^3+4x^4+4x^5+x^7+4x^8+9x^9$ &
  $(1,1,4,7,4,7,8,0,9,0)$ & 4\\
  \hline
  $(5,6,8,4,4,0,0,0,0,0)$ & $6+3x+x^2+3x^3+9x^4+4x^5+10x^6+4x^8+9x^9$ &
  $(1,6,8,2,2,0,8,0,0,0)$ & 4\\
  \hline
  $(1,7,6,8,4,0,0,0,0,0)$ & $7+x+2x^3+4x^5+10x^6+x^7+9x^9$ &
  $(1,9,2,8,4,1,0,0,0,0)$ & 3\\
  \hline
\end{tabular}
\end{narrow}
\end{table}\\ (4). There is received word $u$ whose Lagrange interpolation
polynomial $u(x)=\sum_{i=0}^{q-2}u_i x^i$ satisfies that $u_{q-2}
\neq 0$ and $u_{i}\ne 0$ for all $k \le i \le q-3$. In fact, let
$u=(9, 10, 4, 9, 10, 0, 0, 0,\\ 0, 0)$. Then the Lagrange interpolation
polynomial of $u$ is $u(x)=2+8x+5x^2+7x^3+6x^4+7x^5+x^6+x^7+7x^8+9x^9$.
By computervsearching, we find a codeword $v=(1, 10, 7, 9, 5, 0, 0, 0, 0,
5)$ satisfying $d(u,v)=4 < n-k$. It follows that $u$ is not a deep hole.\\

According to the above examples and exhausting computer search, we propose the
following conjecture as the conclusion of this paper.\\
\\
{\bf Conjecture 4.3.} {\it Assume that the characteristic $p$ of ${\bf F}_q$
is an odd prime. If the Lagrange interpolation polynomial
of the received word $u$ be of the form neither $ax^{q-2}+f_{\le
k-1}(x)$ nor $ax^{k}+f_{\le k-1}(x)$, then $u$ is not a deep hole of
the standard Reed-Solomon code $\mathcal {C}_q({\bf F}_q^*,k)$.} \\
\\
{\bf Acknowledgement.} The authors would like to thank the
anonymous referee for careful reading of the manuscript and
for helpful comments and suggestions. \\

\end{document}